\documentclass[reqno,twoside,a4paper,12pt]{article}
\usepackage[cp866]{inputenc}
\usepackage[english]{babel}
\usepackage {longtable}
\usepackage{amsmath,amssymb,amsfonts}
\usepackage{amsthm}
\usepackage[mathscr]{eucal}
\def\{{\protect\lbrace}
\def\}{\protect\rbrace}

\begin{document}

\begin{center}
\textbf{\Large Semidistributive Laurent Series Rings}

Askar Tuganbaev
\end{center}
\hfill National Research University "MPEI"

\hfill e-mail: tuganbaev@gmail.com

\textbf{Abstract.} If $A$ is a ring with automorphism $\varphi$ and the skew Laurent series ring $A((x,\varphi ))$ is a right semidistributive semilocal ring then $A$ is a right semidistributive right Artinian ring.  The Laurent series ring $A((x))$ is a right semidistributive semilocal ring if and only if $A$ is a right semidistributive right Artinian ring.

The study is supported by Russian Scientific Foundation (project
16-11-10013).

\textbf{Kew words.} Laurent series ring, right semidistributive ring, semilocal ring


\section{Introduction}

All rings are assumed to be associative and with non-zero identity element; all modules are unitary and, unless otherwise specified, all modules are right modules.  The words of type an "Artinian ring" mean "a right and left Artinian ring".

\textbf{1.1. Semidistributive and serial modules and rings.}
A module $M$ is said to be \textsf{distributive} if $X\cap(Y+Z)=X\cap Y+X\cap Z$ for any three its submodules $X,Y,Z$. A module is said to be \textsf{uniserial} if any two its submodules are comparable with respect to inclusion. It is clear that any uniserial module is distributive. The ring $\mathbb{Z}$ of integers is a distributive non-uniserial $\mathbb{Z}$-module. Direct sums of distributive (resp. uniserial) modules are called \textsf{semidistributive} (resp. \textsf{serial}) modules. 

\textbf{1.2. The Laurent series rings $A((x,\varphi ))$ and their modules.}
If $A$ is a ring with automorphism $\varphi$, then
$A((x,\varphi))$ denotes the \textsf{skew Laurent series ring} with
coefficient ring $A$; this ring is formed by all series
$f={\sum\limits_{i=k}^{+\infty}}f_ix^i,$ where $x$ is a variable, $k$ is an
integer (maybe, negative), and all the coefficients $f_i$ are contained in
the ring $A$.  In the ring $A((x,\varphi))$, addition is naturally defined
and multiplication is defined with regard to the relation $xa =
\varphi(a)x$ (for all elements $a\in A$). For $\varphi = 1_A$, we obtain
the ordinary \textsf{Laurent series ring} $A((x))$.\\
For every right $A$-module $M$, we denote by $M((x,\varphi))$ the set of all formal series $\sum _{i =
t}^{\infty}m_ix^i$, where $m_i\in M$, $t\in \mathbb{Z}$, and either $m_t\ne 0$ or $m_i = 0$ for all $i$. The set $M((x,\varphi ))$ is a natural right $A((x,\varphi ))$-module, where module addition is defined
naturally and multiplication by elements of $A((x,\varphi ))$ is defined by the relation 
$$
(\sum _{i = t}^{\infty }m_ix^i)(\sum _{j = s}^{\infty }a_jx^j)= \sum _{k = t+s}^{\infty }(\sum _{i+j = k}m_i\varphi ^i(a_j))x^k.
$$
\textbf{1.3. Remark.} In \cite[Theorem 13.7]{Tug05}, it is proved that $A((x,\varphi ))$ is a right distributive semilocal ring if and only if $A$ is a finite direct product of right uniserial right Artinian rings $A_i$ and $\varphi (A_i)=A_i$ for all $i$. In \cite{Tug20}, it is proved that the ring $A((x,\varphi ))$ is right serial of and only if $A$ is a right serial right Artinian ring. In the both cases, the ring $A((x,\varphi ))$ is a right Artinian ring.

In connection to Remark 1.3, we prove Theorem 1.4 which is the main result of the given paper.

\textbf{1.4. Theorem.} Let $A$ be a ring with automorphism $\varphi$.

\textbf{1.} If $A((x,\varphi ))$ is a right semidistributive semilocal ring, then $A$ is a right semidistributive right Artinian ring and $A((x,\varphi ))$ is a right Artinian ring.

\textbf{2.} Assume that $\varphi(e)=e$ for every local idempotent $e\in A$. Then $A((x,\varphi ))$ is a right semidistributive semilocal ring if and only if $A$ is a right semidistributive right Artinian ring. In this case, $A((x,\varphi ))$ is a right Artinian ring.

\textbf{3.} $A((x))$ is a right semidistributive semilocal ring if and only if $A$ is a right semidistributive right Artinian ring. In this case, $A((x))$ is a right Artinian ring.

In connection to Theorem 1.4, we give Remark 1.5 and Remark 1.6. 

\textbf{1.5. Remark.} Let $F$ be a field and let $A$ be the $5$-dimensional $F$-algebra generated by all $3\times 3$ matrices of the form
$\left(\begin{array}{ccc}
f_{11}& f_{12}& f_{13} \\
0& f_{22}& 0 \\
0& 0& f_{33}
\end{array}\right)$, where $f_{ij}\in F$.
It is directly verified that $A$ is a right semidistributive, left serial, Artinian ring which is not right serial. Therefore, $A((x,\varphi ))$ is not a right serial ring. Thus, it follows from Theorem 1.4 that $A((x))$ is a right semidistributive ring which is not right serial.

\textbf{1.6. Remark.} If $A$ is a field, the the formal power series ring $A[[x]]$ is a right distributive local ring which is not right Artinian.

We present some necessary notation and definitions.\\
Let $A$ be a ring. We denote by $J(A)$ the Jacobson radical of $A$. A ring $A$ is said to be \textsf{semilocal} if $A/J(A)$ is a semisimple
Artinian ring. A ring $A$ is said to be \textsf{local} if $A/J(A)$ is a division ring. 

A module $M$ is said to be \textsf{finite-dimensional} if $M$ does not contain an infinite direct sum of non-zero submodules. A module $M$ is said to be \textsf{quotient finite-dimensional} if all factor modules of the module $M$ are finite-dimensional.

\section{Proof of Theorem 1.4}

\textbf{2.1. Remark.} In \cite[Corollary 5.7]{FacH06}, it is proved that any finite direct sum of quotient finite-dimensional modules is a quotient finite-dimensional module; also see \cite[Corollary 5.24]{Fac19}.

\textbf{2.2. Lemma.} Let $R$ be a ring which has only a finite number of non-isomorphic simple modules.

\textbf{1.} If $M$ is a distributive right $R$-module, then $M$ is finite-dimensional.

\textbf{2.} Every distributive right $R$-module is quotient finite-dimensional.

\textbf{3.} If the ring $R$ is right semidistributive, then every cyclic right $R$-module is quotient finite-dimensional.

\textbf{Proof.} \textbf{1.} Assume the contrary. Then $M$ contains a submodule $X=\oplus_{i=1}^{\infty}X_i$, where $X_i$ is a non-zero-cyclic  module, $i=1,2,\ldots$. For each $i$, the module $X_i$ has a submodule $Y_i$ such that $X_i/Y_i$ is simple. We denote by $Y$ the submodule $\oplus_{i=1}^{\infty}Y_i$ of $X$. The distributive module $M/Y$ contains the submodule $X/Y$ which is isomorphic to the infinite direct sum
$\oplus_{i=1}^{\infty}(X_i/Y_i)$ of simple modules $X_i/Y_i$. Since the ring  $A$ has only a finite number of non-isomorphic simple modules, $X_i/Y_i\cong X_j/Y_j$ for some $i\ne j$. Then $X_i/Y_i\oplus X_j/Y_j$ is a distributive module which is the direct sum of isomorphic simple modules. This is impossible, by \cite{Ste74}.

\textbf{2.} Since all homomorphic image of distributive modules are distributive, the assertion follows from \textbf{1}.

\textbf{3.} The assertion follows from \textbf{2} and Remark 2.1.~\hfill$\square$

\textbf{2.3. Lemma.} If $R$ is a right semidistributive semilocal ring, then every cyclic right $R$-module is quotient finite-dimensional.

Since any semilocal ring has only a finite number of non-isomorphic simple modules, Lemma 2.3 follows from Lemma 2.2(3).

\textbf{2.4. Lemma.} Let $A$ be a ring with automorphism $\varphi$ and $R=A((x,\varphi ))$ the skew Laurent series ring.

\textbf{1.} The ring $A$ is right Artinian if and only if the ring $R$ is right Artinian.

\textbf{2.} If every cyclic right $R$-module is quotient finite-dimensional, then the rings $R$ and $A$ are right Noetherian.

\textbf{3.} If every cyclic right $R$-module is quotient finite-dimensional, then the Jacobson radical $J(R)$ of $R$ is nilpotent.

\textbf{4.} If $R$ is a semilocal ring and every cyclic right $R$-module is quotient finite-dimensional, then the rings $R$ and $A$ are right Artinian.

\textbf{5.} If $R$ is a right semidistributive semilocal ring, then the rings $R$ and $A$ are right Artinian.

\textbf{6.} $A$ is a right Artinian right uniserial ring if and only if $R$ is a right Artinian right uniserial ring, if and only if $R$ is a right uniserial ring.

\textbf{7.} If $A$ is a right Artinian right uniserial ring and $M=mA$ is a cyclic right $A$-module, then the  right $R$-module $mR$ of skew Laurent series is a cyclic uniserial Artinian module.

\textbf{Proof.} \textbf{1.} The assertion is proved in \cite[Proposition 9.2]{Tug05}.

\textbf{2.} The assertion is proved in \cite[Proposition 13.5]{Tug05}.

\textbf{3.} By \textbf{2}, the ring $R$ is right Noetherian. In this case,  the Jacobson radical $J(R)$ is nilpotent, by \cite[Theorem 1(1)]{Tug08}.

\textbf{4.} By \textbf{2} and \textbf{3}, $R$ is a right Noetherian ring with nilpotent Jacobson radical. In addition, $R$ is a semilocal ring, by assumption. It is directly verified that any right Noetherian semilocal ring with nilpotent Jacobson radical is right Artinian. Since the ring $R$ is right Artinian, the ring $A$ is right Artinan, by \textbf{1}.

\textbf{5.} By Lemma 2.3, every cyclic $R$-module is quotient finite-dimensional. Thus the assertion follows from \textbf{4}.

\textbf{6.} The assertion is proved in \cite[Proposition 12.4]{Tug05}.

\textbf{7.} By \textbf{6} $R$ is a right Artinian right uniserial ring. Therefore, $mR\cong R_R/S$, where $S$ is a right ideal of $R$. Therefore, $mR$ is a cyclic uniserial Artinian right $R$-module.~\hfill$\square$

\textbf{2.5. Local idempotents and modules, semiperfect and local rings.} A ring is said to be \textsf{local} if all its non-invertible elements are contained in the Jacobson radical of the ring. For a ring $A$, a right $A$-module $M$ is said to be \textsf{local} if $M$ is a cyclic module and its quotient module modulo its Jacobson radical is simple. For a ring $A$, a non-zero idempotent $e\in A$  is said to be \textsf{local} if $eAe$ is a local ring (equivalently, $eA$ is a local module). A ring $A$ is said to be \textsf{semiperfect} if for its identity element $1_A$, there is a decomposition $1_A=e_1+\cdots+e_n$ into a sum of some orthogonal local idempotents $e_1,\ldots,e_n\in A$; this decomposition is called a \textsf{local decomposition} for the ring $A$. 

In the following familiar assertions \textbf{1}-\textbf{4}, we fix a semiperfect ring $A$ with local decomposition $1_A=e_1+\cdots+e_n$ .

\textbf{1.} If $1_A=f_1+\cdots +f_m$ is one more local decomposition for $A$, then $m=n$ and there is a permutation $\tau$ of the set $\{1,\ldots,n\}$ such that the ring $e_iAe_i$ is isomorphic to the ring $f_{\tau(i)}Af_{\tau(i)}$ and there is an isomorphism of right $A$-modules $e_iA\cong f_{\tau(i)}A$. 

\textbf{2.} If $e$ is a non-zero idempotent of $A$, then there is a non-empty subset $K$ of $\{1,\ldots, n\}$ such that $eA\cong \oplus_{k\in K}e_kA$.

\textbf{3.} A right $A$-module $M$ is distributive if and only if $Me_i$ is a uniserial right $e_iAe_i$-module, for each $e_i$.

\textbf{4.} The ring $A$ is right semidistributive if and only if $e_jAe_i$ is a uniserial right $e_iAe_i$-module, for each $e_i$ and $e_j$.

\textbf{5.} The ring $A$ is right semidistributive if and only if the right $A$-module $e_iA$ is distributive for each $i$, if and only if for any local decomposition $1_A=f_1+\cdots +f_m$, the right $A$-module $f_iA$ is distributive for each $i$, if and only if for any local decomposition $1_A=f_1+\cdots +f_m$, the right $f_iAf_i$-module $f_jAf_i$ is uniserial for each $i$.

\textbf{Proof.} \textbf{1, 2.} The assertions are well known; e.g., see \cite[Section 27]{AndF92} or \cite[Section 6.3]{Tug98}.

\textbf{3.} The assertion is proved in \cite[Lemma 4]{Ful78}.

\textbf{4.} The assertion follows from \textbf{3}.

\textbf{5.} The assertion follows from \textbf{1} and \textbf{4}.~\hfill$\square$

\textbf{2.6. Lemma.} Let $A$ be a ring, $\varphi$ be an automorphism of $A$ such that $\varphi(e)=e$ for every idempotent $e\in A$, and let $R=A((x,\varphi ))$ be the skew Laurent series ring.

\textbf{1.} For any non-zero idempotent $e\in A$ and each right $A$-module $M$, the skew Laurent series ring $(eAe)((x,\varphi))$ is naturally isomorphic to the ring $eRe$ and the right $(eAe)((x,\varphi))$-module $(Me)((x,\varphi))$ of skew Laurent series can be naturally identified with the right $eRe$-module $(Me)((x,\varphi))$.

\textbf{2.} If $e$ is a non-zero idempotent of the ring $A$ such that $eAe$ is a right uniserial right Artinian ring, then the ring $eRe$ is a right uniserial right Artinian ring and $e$ is a local idempotent of $R$.

\textbf{3.} If $A$ is a right semidistributive right Artinian ring, then $R$ is a right semidistributive right Artinian ring.

\textbf{Proof.} \textbf{1.} The assertion is directly verified.

\textbf{2.} By 2.4(6), the skew Laurent series ring $(eAe)((x,\varphi))$ is a right uniserial right Artinian ring. By \textbf{1}, the ring $eRe$ is a right uniserial right Artinian ring. Therefore, $e$ is a local idempotent of $R$.

\textbf{3.} Since $A$ is a right Artinian ring, it follows from Lemma 2.4(1) that $R$ is a right Artinian ring. In particular, the ring $A$ is semiperfect and its identity element $1_A$ has a decomposition $1_A=e_1+\cdots+e_n$ into a sum of some orthogonal local idempotents $e_1,\ldots,e_n\in A$. Since $A$ is a right semidistributive semiperfect ring, it follows from 2.5(4) that $e_jAe_i$ is a uniserial right $e_iAe_i$-module, for each $e_i$ and $e_j$. Since $A$ is a right Artinian ring, it is directly verified that each ring $e_iAe_i$ is right Artinian.  By \textbf{2},  each $e_i$ is a local idempotent of $R$. 
By \textbf{1} and Lemma 2.4(6), the skew Laurent series ring $(e_iAe_i)((x,\varphi))$ is naturally isomorphic to the ring $e_iRe_i$ and is right uniserial.  Since $e_iAe_i$ is a right uniserial right Artinian ring, all cyclic right $e_iRe_i$-modules are uiserial Artinian right modules, by Lemma 2.4(7). By \textbf{1}, the right $(e_iAe_i)((x,\varphi))$-module $(e_jAe_i)((x,\varphi))$ of skew Laurent series can be naturally identified with the right $e_iRe_i$-module $(e_jAe_i)((x,\varphi))$.
Since $1_R=1_A=e_1+\cdots+e_n$ is the sum of orthogonal local idempotents $e_1,\ldots,e_n\in R$, it follows from 2.5(4) that $R$ is a right semidistributive ring.~\hfill$\square$

\textbf{2.7. The completion of the proof of Theorem 1.4.} Let $R=A((x,\varphi ))$.

\textbf{1.} Let $R$ be a right semidistributive semilocal ring. By Lemma 2.4(5), the rings $R$ and $A$ are right Artinian.

Let $\{e_1,\ldots, e_n\}$ be a complete set of local orthogonal idempotents of the right Artinian right semidistributive ring $A$. By 2.5(4), each of the rings $e_iAe_i$ are right Artinian right uniserial rings. By Lemma 2.4(6), each of the rings $e_iRe_i$ are right Artinian right uniserial rings. In parfticular, each of the rings $e_iRe_i$ are right Artinian right uniserial rings and $\{e_1,\ldots, e_n\}$ is a complete set of local orthogonal idempotents of the right Artinian ring $R$. 

By applying 2.5(5) to the right semidistributive right Artinian ring $R$, we obtain that  for all $i$ and $j$, the right $e_iRe_i$-module $e_jRe_i$ is uniserial. We fix $i$ and $j$. By 2.5(5), it is sufficient to prove the right $e_iAe_i$-module $e_jAe_i$ is uniserial. 

Let $e_jae_i,e_jbe_i\in e_jAe_i$ and $e_jae_i\notin e_jbe_iA$. 
If $e_jae_i\in e_jbe_iR$, then $e_jae_i=e_jbe_if$ for some $f\in R$. Then $e_jae_i=e_jbe_if_0$ for the constant term $f_0$ of $f$ and $e_jae_i\in e_jbe_iA$; this is a contradiction. Since the right $e_iRe_i$-module $e_jRe_i$ is uniserial, we have $e_jbe_i\in e_jae_iR$ and $e_jbe_i=e_jae_ig$ for some $g\in R$. Then Then $e_jbe_i=e_jae_ig_0$ for the constant term $g_0$ of $g$ and $e_jbe_i\in e_jae_iA$. Therefore, the right $e_iAe_i$-module $e_jAe_i$ is uniserial. By 2.5(5), $A$ is a right semidistributive right Artinian ring.

\textbf{2.} The assertion follows from \textbf{1} and Lemma 2.6(2).

\textbf{3.} Let $A$ be a right semidistributive right Artinian ring. By Lemma 2.6(3), $R$ is a right semidistributive right Artinian ring.

Let $R$ be a right semidistributive right Artinian ring. By \textbf{1}, $A$ is a right semidistributive right Artinian ring.~\hfill$\square$

\section{Open Questions}

Let $A$ be a ring with automorphism $\varphi$.

\textbf{3.1.} Let $A$ be a right semidistributive right Artinian ring. Is it true that $A((x,\varphi ))$ is a right semidistributive ring?

\textbf{3.2.} Let $R=A((x,\varphi ))$ be a \textsf{regular} ring, i.e., $r\in rRr$ for each $r\in R$. Is it true that the ring $R$ is Artinian? This is true if the  automorphism $\varphi$ is of finite order; see \cite[Theorem 1]{Son95}.

\textbf{3.3.} Let $A$ be a ring such that the ring $A((x,\varphi ))$ is semilocal. Is it true that $A$ is semiperfect and the Jacobson radical of $A$ is nil? This is true if $\varphi =1_A$; see \cite{Zie14}.

\textbf{3.4.} When is the ring $A((x,\varphi ))$ semilocal?

\textbf{3.5.} When is the ring $A((x,\varphi ))$ right distributive?

The author thanks Alberto Facchini for Remark 2.1.

\end{document}